\documentclass[a4,11pt]{article}
\usepackage[centertags]{amsmath}
\usepackage{amsmath,amsfonts,amssymb,amsthm}
\usepackage[latin1]{inputenc}
\newcommand{\bb}{\begin{equation}}
\newcommand{\ee}{\end{equation}}
\newcommand{\ds}{\displaystyle }

\newcommand{\lh}{\Delta_{H^{1}}}
\def\p{\partial}

\def\x1{{\xi }_{xx}}
\def\x2{{\xi }_{yy}}
\def\x3{{\xi }_{xy}}

\def\e1{{\eta }_{xx}}
\def\e2{{\eta }_{yy}}
\def\e3{{\eta }_{xy}}

\def\al{\alpha}

\def\var{\varphi}

\def\de{\delta}

\def\be{\beta}
\def\De{\Delta}

\def\l1{{\lambda}_1}

\def\kd{\partial}

\begin{document}
\title{\huge \bf Group Classification of Semilinear Kohn-Laplace
Equations}
\author{\rm \large Yuri Bozhkov and Igor Leite Freire\\ \\
\it Instituto de Matem\'atica,
Estat\'\i stica e \\ \it Computa\c c\~ao Cient\'\i fica - IMECC \\
\it Universidade Estadual de Campinas - UNICAMP \\ \it C.P.
$6065$, $13083$-$970$ - Campinas - SP, Brasil
\\ \rm E-mail: bozhkov@ime.unicamp.br \\ \ \ \ \ \
igor@ime.unicamp.br }
\date{\ }
\maketitle
\vspace{1cm}
\begin{abstract}
We study the Lie point symmetries of semilinear Kohn-Laplace
equations on the Heisenberg group $H^1$ and obtain a complete
group classification of these equations.
\end{abstract}


\section{Introduction}

\

 The Heisenberg group $H^n$ topologically is the real vector space
 ${\mathbb{R} }^{2n+1}$.
 Its Lie group structure is determined by the product
 \[(x,y,t)  (x^0,y^0,t^0) = (x+x^0, y+y^0, t+t^0 +2\sum_{i=1}^{n}(y_i x^0_i-x_i y^0_i)), \]
 where $(x,y,t),(x^0,y^0,t^0)\in {\mathbb{R} }^n\times {\mathbb{R} }^n\times \mathbb{R} =H^n$.
 It is easy to verify that the operators
 \[ T=\frac{\kd }{\kd t},\;\;X_i=\frac{\p}{\p x_i}+2y \frac{\p}{\p
 t},\;\;Y_i=\frac{\p}{\p y_i}-2x \frac{\p}{\p t},\]
 where $i=1,2,...,n,$ form a basis of the left-invariant vector fields on $H^n$ and
 satisfy the following commutation relations:
 \[ [X_i,Y_i]=-4 {\de }_{ij} T,\;\;\; [X_i,X_j]=[Y_i,Y_j]=[X_i,T]=[Y_i,T]=0. \]
 These formulae present in an abstract form the commutation relations for the quantum-mechanical position
 and momentum operators in $n-$dimensional configuration space.
 This justifies the name Heisenberg group.

 In the last few decades the Heisenberg group $H^n$ was intensively and extensively
 studied by a considerable number of authors using methods and
 approaches which come from algebraick and differential geometry, real and complex analysis,
 mathematical physics and applications. A big part of the corresponding works treats
 partial differential equations on $H^n$. In this regard various authors have obtained existence
 and nonexistence results for equations involving Kohn-Laplace operators. Recall that the Kohn-Laplace
 operator ${\De }_{H^n}$ is the natural subeliptic Laplacian on $H^n$ defined by
 \[ {\De }_{H^n}=\sum_{i=1}^{n} (X^2_i+Y^2_i). \]
 Although there are similarities between $ {\De }_{H^n}$ and the classical Laplacian they are essentially
 different. E.g. the Kohn-Laplace operator is not a strongly elliptic operator. It is a typical representative
 of the hypoeliptic operators (\cite{lars}). (Since the study of
 hypoellipticity properties is not subject of this paper we shall
 not comment more on this point.)

 In \cite{gl3} Garofalo and Lanconelli established existence, regularity and nonexistence results for the
 Kohn-Laplace equation
 \[ {\De }_{H^n}u+f(u)=0\]
 in an open bounded or unbounded subset of $H^n$ with homogeneous Dirichlet boundary condition. One of the
 motivations to study such semilinear equations is the fact that they may arise as Euler-Lagrange
 equations in some variational problems on Cauchy-Riemann (CR) manifolds as in the works of Jerison and Lee
 \cite{jl1,jl2} on the CR Yamabe problem.
 The existence of weak solutions is proved in \cite{gl3} provided the nonlinear term satisfies some growth
 conditions of the form $f(u)= o(|u|^{(Q+2)/(Q-2)})$ as $|u|\rightarrow\infty $,
 where $Q=2n+2$ is the so-called homogeneous dimension of
 $H^n$ (\cite{f2}). The exponent $(Q+2)/(Q-2)$ is the critical exponent
 for the Stein's Sobolev space (\cite{jl2}). The nonexistence results
 follow from remarkable Pokhozhaev Identities established in \cite{gl3} for the solutions of Kohn-Laplace equations on
 the Heisenberg group. The Dirichlet problem for the Kohn
Laplacian on $H^n$
 was studied before by Jerison in \cite{j1,j2}. See also \cite{bia} for existence of classical nonnegative solutions of
 semilinear Kohn-Laplace equations. General nonexistence results for solutions of semilinear
differential inequalities on the Heisenberg group were obtained by
Pokhozhaev and Veron in \cite{pv}. Since there is a huge number of
works dedicated to Heisenberg groups (see \cite{ams}) and the
study of PDE on $H^n$, in order not to
 increase the volume of this paper, we shall not present here further details, directing the interested
 reader to the already cited works as well as to \cite{by,be,ib1,ib2,f1,f2,ho,hm} and the references therein.

 The purpose of the present paper is to enlighten the properties of the Kohn-Laplace equations from
 the point of view of the S. Lie Symmetry Theory, which to our knowledge has
 not been previously done. We shall obtain complete group classification of semilinear partial
 differential equations on $H^1$ of the following form
\bb\label{001} \lh u+f(u)=0, \ee
 where $\lh $ is the Kohn-Laplace operator on $H^1$ and $f$ is a generic function.

 The importance of group classification of differential equations was first
 emphasized by Ovsiannikov in 1950s-1960s, when he and his school
 began a systematic research program of successfully applying
 modern group analysis methods to wide range of physically
 important problems. Following Olver (\cite{ol}, p. 182), we recall that to perform a group
classification on a differential equation
 involving a generic function $f$ consists of finding the Lie point symmetries of the given equation with
 arbitrary $f$, and, then, to determine all possible particular forms of $f$ for which the symmetry group
 can be enlarged. It is worth observing that for problems which arise from physics, quite often there
 exists a physical motivation for considering such specific cases.

 The Heisenberg group
 $H^1$ itself possesses the rich properties of $H^n$ (see \cite{be}) and the calculations of the symmetry group of this
 model problem give insights for the general case $n>1$. For this reason, and for the sake of simplicity
 and clarity we restrict ourselves to $H^1$.

 We write the Kohn-Laplace operator
 \bb \label{m01} \lh =X^2+Y^2, \ee
 where
 \bb\label{004} X=\frac{\p}{\p x}+2y \frac{\p}{\p t} \ee
 and
 \bb\label{005} Y=\frac{\p}{\p y}-2x \frac{\p}{\p t} .\ee
 Then the equation (\ref{001}) for $u=u(x,y,t)$ in more details
 reads
 \bb\label{002}
 u_{xx}+u_{yy}+4(x^{2}+y^{2})u_{tt}+4yu_{xt}-4xu_{yt}+f(u)=0 .\ee

 We shall not present preliminaries concerning  Lie point symmetries of
 differential equations supposing that the reader is familiar with
 the basic notions and methods of contemporary group analysis \cite{bk,ib,ol,ov}.

 The main result in this paper is the following

\textbf{Theorem} \it The widest Lie point symmetry group of the
Kohn-Laplace equation $(\ref{001})$ with an arbitrary $f(u)$ is
determined by the operators
\bb\label{u01} T=\frac{\kd }{\kd t},\;\;\; R= y\frac{\kd }{\kd
x}-x \frac{\kd }{\kd y},\;\;\; \tilde{X}=\frac{\p}{\p x}-2y \frac{\p}{\p
t}, \;\;\;\tilde{Y}=\frac{\p}{\p y}+2x \frac{\p}{\p t}, \ee
that is, by a translation in t, a rotation in the x-y plane and the generators
of right multiplication in the Heisenberg group $H^1$.

 For some special choices of the right-hand side $f(u)$ it can be extended in the
 cases listed below. We shall write only
 the generators additional to $(\ref{u01})$.

 $(i)$ If $f(u)=0$, then
 \bb
 \label{u04} V_1 = (xt-x^{2}y-y^{3})\frac{\p }{\p x} +
 (yt+x^{3}+xy^{2})\frac{\p }{\p y} + (t^{2}-(x^{2}+y^{2})^{2})\frac{\p }{\p t}-
 t u \frac{\p }{\p u}\ee
  \bb\label{u05} V_2  = (t-4xy)\frac{\p }{\p x} +
 (3x^{2}-y^{2})\frac{\p }{\p y} - (2yt+2x^{3}+2xy^{2})\frac{\p }{\p t}
 + 2 y u \frac{\p }{\p u} \ee
  \bb \label{u06} V_3 = (x^{2}-3y^{2})\frac{\p }{\p x} +
 (t+4xy)\frac{\p }{\p y} + (2xt-2x^{2}y-2y^{3})\frac{\p }{\p t}
 - 2x u \frac{\p }{\p u}\ee
  \bb \label{u02} Z_1 =x\frac{\p}{\p x}+y\frac{\p}{\p y}+2t
\frac{\p}{\p t} ,  \;\;\; Z_2 = u\frac{\kd }{\kd u},\;\;\;W= \be
(x,y,t)\frac{\kd }{\kd u},
 \ee
 where $\lh \be =0$.

 $(ii)$ If $f(u)=c=const $, then this case is reduced to $(i)$
 by the change $u=v+ c x^2/2$.

 $(iii)$ If $f(u)=k.u $, k-constant, then
 \bb \label{u08} Z_2 = u\frac{\kd }{\kd u},\;\;\; W= \be (x,y,t)\frac{\kd }{\kd u}, \ee
 where $\lh \be +k\be=0$.

 $(iv)$ If $f(u)=k.u^{\ds p} $, $p\neq 0, p\neq
 1$, we have the generator of dilations
 \bb\label{042} Z= x\frac{\p}{\p x}+y\frac{\p}{\p y}+2t
 \frac{\p}{\p t}+\frac{2}{1-p}u\frac{\p}{\p u}. \ee In the critical
 case $f(u)=k.u^{\ds 3} $, there are three additional generators, namely
 $V_1,V_2,V_3$ given in $(\ref{u04}), (\ref{u05}),(\ref{u06})$
 respectively.

 $(v)$ If $f(u)=k.e^{\ds u}$ then the operator
 \bb\label{040} Z_3=x\frac{\p}{\p x}+y\frac{\p}{\p y}+2t
\frac{\p}{\p t}-2\frac{\p}{\p u} \ee
 generates a sub-group of the Lie point symmetry group of
 $(\ref{001})$.

\rm

\

 This classification is similar to that for semilinear equations in ${\mathbb{R} }^n$ involving
 Laplace or polyharmonic operators \cite{sv}. We also observe that for power nonlinearity $f(u)=ku^p$
 exactly in the critical case $p=3= (Q+2)/(Q-2)$, $Q=2.1 +2=4$ being the homogeneous dimension of $H^1$,
 the symmetry group is expanded by three aditional generators (see $(iv)$ of the main theorem). This fact
 suggests that in the critical case maybe there are further properties as pointed out in \cite{yb1,yb2}
 regarding other differential equations. This is our motivation to use the above group classification in
 three subsequent papers \cite{yi2,yi3,yi4}. In \cite{yi2} we study the variational properties of
 Kohn-Laplace equations and we find out which of the already found Lie point
 symmetries are variational/divergence symmetries. Further in \cite{yi4} we establish the
 corresponding conservation laws via the Noether Theorem. In \cite{yi3} we discuss the
 invariant solutions of various Kohn-Laplace
 equations on the Heisenberg group.

 The group classification of Kohn-Laplace equations on the Heisenberg group $H^n$, $n>1$, will be treated
 elswhere.

 This paper is organized as follows. In the next section we obtain the determining equations for the Lie
 point symmetries of the equation (\ref{001}). This process is essentially simplified by the use of two
 theorems of Bluman \cite{b1,bk}. Then in section 3 we obtain some formulae which are consequences of the
 determining equations. They are used in the proof of the main theorem, given in sections 4-9.

 \section{The determining equations}

 In this section we obtain the determining equations for a Lie point symmetry of the
 Kohn-Laplace equation (\ref{001}) with infinitesimal generator
 \bb\label{006} S=\xi\frac{\p}{\p x}+\phi \frac{\p}{\p y}
+\tau\frac{\p}{\p t}+\eta \frac{\p}{\p u}.\ee

 To begin with, we observe that the symmetry calculation is drastically simplified if we apply
 two theorems of Bluman \cite{b1,bk}. Indeed, Theorem 4.2.3-1, \cite{bk}, p. 174, implies that
 $\xi , \phi $ and $\tau $ do not depend on $u$. Then by Theorem 4.2.3-6, \cite{bk}, p. 175, we
 conclude that $\eta $ is a linear function of $u$. Therefore the infinitesimals are of the
 following form
 \bb\label{003} \left\{\begin{array}{r c l} \xi & = & \xi(x,y,t),\\
 \phi & = & \phi(x,y,t),\\
 \tau & = & \tau(x,y,t),\\
 \eta & = & \al(x,y,t)u+\be(x,y,t), \end{array}\right. \ee
 where $\al =\al (x,y,t)$ and $\be =\be(x,y,t)$ are functions to be determined.

 We denote
 \[ H:=\lh u + f(u). \]
 The equation (\ref{001}) admits the symmetry (\ref{006}) if and only if
 \[ \hat{S}H=0 \]
 when $H=0$ (\cite{bk,ol}), where
\bb\label{007} \begin{array}{lll} \hat{S}&=&\xi\ds{\frac{\p}{\p
x}}+\phi \ds{\frac{\p}{\p y}} +\tau\frac{\p}{\p t}+\eta
\ds{\frac{\p}{\p u}}+\eta^{(1)}_{x}\frac{\p}{\p
u_{x}}\\
& & \\
&+&\eta^{(1)}_{y}\ds{\frac{\p}{\p u_{y}}}+\eta_{t}^{(1)}\frac{\p
}{\p u_{t}}+\eta^{(2)}_{xx}\frac{\p}{\p
u_{xx}}+\eta_{yy}^{(2)}\frac{\p}{\p
u_{yy}}\\
& & \\
&+&\eta_{xt}^{(2)}\ds{\frac{\p}{\p
u_{xt}}}+\eta_{yt}^{(2)}\ds{\frac{\p}{\p u_{yt}}}+
\eta_{tt}^{(2)}\ds{\frac{\p}{\p u_{tt}}}
+\eta_{xy}^{(2)}\ds{\frac{\p}{\p u_{xy}}} \end{array}\ee
 is the second order extension of $S$ (\cite{bk,ol}). Then the symmetry condition can be
 written as
 \bb\label{011}\begin{array}{lll}
(8x\xi&+&8y\phi)u_{tt}+4\phi u_{xt}-4\xi u_{yt}+\eta
f^{'}(u)\\ \\
&+&\eta^{(2)}_{xx}+\eta_{yy}^{(2)}+4(x^{2}+ y^{2})\eta_{tt}^{(2)}
+ 4y\eta^{(2)}_{xt}-4x\eta_{yt}^{(2)}=0,\end{array} \ee
 when $H=0$. (The subscripts denote partial derivatives, e.g. $u_x=\frac{\p u}{\p x}$.
 Only in the extension coefficients like $\eta^{(j)}_{xt}$ the subscripts mean indices. We also
 suppose that the considered functions are sufficiently smooth in order that the derivatives we write
 to exist.) Further, using the corresponding formulae for the extended infinitesimals (\cite{bk,ol})
 we calculate
 \bb\label{008}
\eta^{(1)}_{x}=\be_{x}+\al_{x}u+(\al-\xi_{x})u_{x}-\phi_{x}u_{y}-\tau_{x}u_{t},
\ee \bb\label{009}
\eta^{(1)}_{y}=\be_{y}+\al_{y}u-\xi_{y}u_{x}+(\al-\phi_{y})u_{y}-\tau_{y}u_{t},
\ee \bb\label{010}
\eta^{(1)}_{t}=\be_{t}+\al_{t}u-\xi_{t}u_{x}-\phi_{t}u_{y}+(\al-\tau_{t})u_{t}.
\ee \bb\label{012}\begin{array}{lll}
\eta_{xx}^{(2)}& = & \be_{xx}+\al_{xx}u+(2\al_{x}-\xi_{xx})u_{x}-\phi_{xx}u_{y}-\tau_{xx}u_{t}\\
& + &
(\al-2\xi_{x})u_{xx}-2\phi_{x}u_{xy}-2\tau_{x}u_{xt},\end{array}
\ee \bb\label{013}\begin{array}{lll}
\eta_{yy}^{(2)}&=&\be_{yy}+\al_{yy}u-\xi_{yy}u_{x}+(2\al_{y}-\phi_{yy})u_{y}
-\tau_{yy}u_{t}\\
&-&2\xi_{y}u_{xy}+(\al-2\phi_{y})u_{yy}-2\tau_{y}u_{yt},\end{array}
\ee \bb\label{014}\begin{array}{lll}
\eta_{tt}^{(2)}&=&\be_{tt}+\al_{tt}u-\xi_{tt}u_{x}-\phi_{tt}u_{y}
+(2\al_{t}-\tau_{tt})u_{t}\\&-&2\xi_{t}u_{xt}-2\phi_{t}u_{yt}+(\al-2\tau_{t})u_{tt},\end{array}
\ee \bb\label{015}\begin{array}{lll}
\eta_{xt}^{(2)}&=&\be_{xt}+\al_{xt}u+(\al_{t}-\xi_{xt})u_{x}-\phi_{xt}u_{y}+(\al_{x}-\tau_{xt})u_{t}\\
&-&\xi_{t}u_{xx}-\phi_{t}u_{xy}-\phi_{x}u_{yt}+(\al-\xi_{x}-\tau_{t})u_{xt}-\tau_{x}u_{tt},\end{array}
\ee \bb\label{016}\begin{array}{lll}
\eta_{yt}^{(2)}&=&\be_{yt}+\al_{yt}u-\xi_{yt}u_{x}
+(\al_{t}-\phi_{yt})u_{y}
+(\al_{y}-\tau_{yt})u_{t}\\
&-&\xi_{y}u_{xt}
-\xi_{t}u_{xy}-\phi_{t}u_{yy}+(\al-\phi_{y}-\tau_{t})u_{yt}-\tau_{y}u_{tt}.\end{array}
\ee

Now substituting $(\ref{008})-(\ref{016})$ into the symmetry
condition (\ref{011}), after some tedious work, we obtain
\bb\label{017}\begin{array}{l}\displaystyle{(\al
u+\be)f^{'}(u)+\lh
\be +(\lh \al)u+[2\al_{x}+4y{\al }_{t}-\lh \xi]u_{x}}\\ \\
\displaystyle{+[2\al_{y}+4x{\al }_{t}-\lh
\phi]u_{y}+[8(x^{2}+y^{2})\al_{t}+4y\al_{x}-4x\al_{y}-\lh
\tau]u_{t}}\\ \\
\displaystyle{+[-2\phi_{x}-2\xi_{y}-4y\phi_{t}+4x\xi_{t}]u_{xy}+[\al-2\xi_{x}-4y\xi_{t}]u_{xx}+
[\al-2\phi_{y}+4x\phi_{t}]u_{yy}}\\ \\
\displaystyle{ +[8x\xi + 8y\phi
+4(x^{2}+y^{2})(\al-2\tau_{t})-4y\tau_{x}+4x\tau_{y}]u_{tt}}\\ \\
\displaystyle{+[4\phi-2\tau_{x}-8(x^{2}+y^{2})\xi_{t}+4y(\al-\xi_{x}-\tau_{t})+4x\xi_{y}]u_{xt}}\\ \\
\displaystyle{+[-4\xi-2\tau_{y}-8(x^{2}+y^{2})\phi_{t}-4y\phi_{x}-4x(\al-\phi_{y}-\tau_{t})]u_{yt}=0,}
\end{array}\ee
when $H=0$. Then, expressing $u_{xx}$ from (\ref{002}) and
substituting in (\ref{017}), we obtain an identity for all values
of $(x,y,t,u, u_x,u_y,u_t,u_{xy},u_{xt},u_{yy},u_{yt},u_{tt})$.
Equating to zero the coefficients of the derivatives of $u$ and
the free term, we obtain the following nine determining equations:
\bb\label{018} \xi_{x}+2y\xi_{t}-\phi_{y}+2x\phi_{t}=0, \ee
\bb\label{019} \xi_{y}-2x\xi_{t}+\phi_{x}+2y\phi_{t}=0, \ee
\bb\label{020} \lh \xi=2X\al, \ee \bb\label{021} \lh \phi=2Y\al,
\ee \bb\label{022} \lh \tau=4yX\al-4xY\al, \ee \bb\label{023} \al
u f^{'}(u)+\be f^{'}(u)+(\lh \al)u+\lh
\be+[4y\xi_{t}+2\xi_{x}-\al]f(u)=0, \ee \bb\label{024}
2y\xi_{x}+2x\xi_{y}+4(y^{2}-x^{2})\xi_{t}+2
\phi-\tau_{x}-2y\tau_{t}=0, \ee \bb\label{025}
4x\xi_{x}+8xy\xi_{t}+2\xi+2y\phi_{x}-2x\phi_{y}+4(x^{2}+y^{2})\phi_{t}+\tau_{y}-2x\tau_{t}=0,\ee
\bb\label{026}
2(x^{2}+y^{2})\xi_{x}+4y(x^{2}+y^{2})\xi_{t}+2x\xi+2y\phi-y\tau_{x}+x\tau_{y}-2(x^{2}+y^{2})\tau_{t}=0,
\ee
 where the operators $X$ and $Y$ are defined by (\ref{004}) and (\ref{005}). Multiplying
 correspondingly the equations (\ref{018}), (\ref{019}), (\ref{024}) and (\ref{025}) we obtain a
 relation which symbolically can be written as
 ``(\ref{026})=$y$.(\ref{024})+ $x$.(\ref{025})- $x$.(\ref{018})-$y$.(\ref{019})." Hence,
 the equation (\ref{026}) is a consequence of (\ref{018}), (\ref{019}), (\ref{024}) and
 (\ref{025}). Another straightforward calculation shows that (\ref{022}) also follows from these
 equations. Therefore there are seven independent determining equations which in terms of the
 operators $X$ and $Y$ can be written in the following simplified form:
\bb\label{027} X\xi-Y\phi=0, \ee \bb\label{028} Y\xi+X\phi=0, \ee
\bb\label{029} \lh \xi=2X\al, \ee \bb\label{030} \lh \phi=2Y\al,
\ee \bb\label{031} \al u f^{'}(u)+\be f^{'}(u)+(\lh \al)u+\lh
\be+(2X\xi-\al)f(u)=0, \ee \bb\label{032}
X\tau=2yX\xi+2xY\xi+2\phi, \ee \bb\label{033}
Y\tau=-2xX\xi+2yY\xi-2\xi. \ee

We conclude this section by noting that the system of two
equations (\ref{027}) - (\ref{028}) may be considered as a
Heisenberg group generalization of the Cauchy-Riemann equations.

\section{Some consequences of the determining equations}

 \bf Proposition 1. \it If the infinitesimals $\xi $ and $\phi $
satisfy $(\ref{027})-(\ref{028})$ then \bb\label{u09} \lh \phi = 4
{\xi }_t \ee and \bb\label{u10} \lh \xi = -4 {\phi }_t. \ee

\

\bf Proof. \rm  We apply $X$ to equation (\ref{028}), $Y$ to
equation (\ref{027}) and subtract the resulted equations. In this
way we obtain
\[ (XY-YX)\xi + (X^2+Y^2)\phi =0,\]
 which implies (\ref{u09}) since the commutator
 \bb\label{m03} [X,Y]=-4{\kd }_t \ee
 and $\lh = X^2+Y^2$. The equation (\ref{u10}) can be derived in an analogous way.

 \

 \bf Corollary. \it If $\al , \xi , \phi $ satisfy $(\ref{027})-(\ref{030})$, then
\bb\label{034} X\al=-2\phi_{t}, \ee \bb\label{035} Y\al=2\xi_{t}.
\ee

\

\bf Proposition 2. \it If $\xi , \phi $ and $\tau $ satisfy
$(\ref{027}),(\ref{028}),(\ref{032})$ and $(\ref{033})$, then
\bb\label{036} \tau_{t}=2y\xi_{t}-2x\phi_{t}+2X\xi. \ee

\

\bf Proof. \rm We just sketch the proof. We apply the operator $X$
to equation (\ref{033}), the operator $Y$ to equation (\ref{032})
and subtract. The resulted equation, with the use of the
commutator (\ref{m03}) and (\ref{027}), (\ref{028}), leads to
(\ref{036}).

\

\bf Proposition 3. \bb \label{m04} {\al }_t=-(X\xi )_t . \ee

\

\bf Proof. \rm We apply the operator $X$ to (\ref{035}), the
operator $Y$ to (\ref{034}) and subtract. Then by \bb\label{m06}
[X,{\kd }_t]=[Y,{\kd }_t]=0, \ee (\ref{m01}) and (\ref{028}) we
obtain (\ref{m04}).

\section{The Lie point symmetries for arbitrary $f(u)$}

In this section we prove the main theorem for general right-hand
side of the Kohn-Laplace equation (\ref{001}).

Since $f(u)$ is an arbitrary function, then $\al = \be =0$ by
(\ref{031}). Thus the equations (\ref{029}) and (\ref{030}) imply
that \bb\label{t001} \lh \xi = \lh\phi =0.\ee Then by
(\ref{t001}), (\ref{u09}) and (\ref{u10}) it follows that
$${\xi }_t={\phi }_t=0 $$ and therefore $\xi $ and $\phi $ are
functions of $x$ and $y$ only. On the other hand, from the
equation (\ref{031}), $$2X\xi -\al =0,$$ which implies that $\xi $
depends only on $y$ since $\al =0$ and $\xi $ does not depend on
$t$. From the equation (\ref{027}) for $\xi (y)$ and $\phi (x,y)$
it follows that $\phi $ depends only on $x$. Further, from
(\ref{t001}), we obtain that $$\xi = a_1 y +a_2,$$
 $$\phi = A x + a_3,$$ where $A,a_1, a_2,a_3$ are arbitrary
constants. Now from (\ref{028}) we get that $A=-a_1$, that is
$$\phi = -a_1 x + a_3.$$ Substituting $\xi $ and $\phi $ into
(\ref{032}) and (\ref{033}) gives  $$X\tau - 2 a_3=0,$$ $$Y\tau +
2 a_2=0.$$ Applying $Y$ and $X$ to the latter two equations,
subtracting and using the commutator $[X,Y]=-4{\p }_t$ we obtain
that $\tau $ does not depend on $t$. Therefore $${\tau }_x-
2a_3=0,$$
 $${\tau }_y + 2a_2=0,$$ from which we conclude easily that $\tau = 2a_3 x - 2 a_2 y + a_4$ where $a_4$ is
 another arbitrary constant. In this way
\bb\label{038} \left\{\begin{array}{r c l } \xi & = & a_{1}y
+ a_{2},   \\
\phi & = &- a_{1}x  + a_{3},   \\
\tau & = & 2a_{3}x - 2a_{2}y + a_{4},\\
\eta & = & 0,\end{array} \right.\ee which proves the first
statement of the main theorem.

 \section{The Lie point symmetries for $ f(u)=ke^{\ds{u}}$}

 In this section we prove the item (v) of the main theorem.

 We substitute $f(u)=ke^u$ into (\ref{031}):
 \[ \al k u e^u +\be k e^u + (\lh \al)u + (\lh\be ) +[2X\xi -\al ]ke^u=0.\]
 Hence
 \bb\label{m07} \al =0 \ee
 and
 \bb\label{m08} \be +2X\xi -\al =0. \ee
 From (\ref{m07}) and (\ref{034}), (\ref{035}) it follows that $\xi =\xi (x,y)$ and
 $\phi =\phi (x,y)$. From (\ref{m08}):
 \bb\label{m09} \be +2{\xi }_x=0 \ee
 and hence $\be =\be (x,y)$.

 Further, the relation (\ref{036}) implies that
 \[ {\tau }_t=2 {\xi }_x \]
 since $\xi $ and $\phi $ do not depend  on $t$. Therefore there exists a function $h(x,y)$ such
 that
 \bb\label{m10} \tau =2 t {\xi }_x +h(x,y). \ee
 We observe that he right-hand side of the equation (\ref{032}) does not depend on $t$, while
 the left-hand side is
 \[ 2t\;{\xi }_{xx}+h_x +4y{\xi }_x.\]
 Thus
 \bb\label{m11} {\xi }_{xx}=0.\ee
 Analogously, from (\ref{033}), we obtain that
 \bb\label{m12} {\xi }_{xy}=0.\ee
 On the other hand $\lh\xi =0$ by (\ref{029}) since $\al =0$ (see (\ref{m07})). But $\xi =\xi
 (x,y)$ and therefore $\xi $ is harmonic:
 \bb\label{m13} {\xi }_{xx}+{\xi }_{yy}=0. \ee
 The equations (\ref{m11}), (\ref{m12}) and (\ref{m13}) imply that
 \[ \xi = a_5x + a_1 y + a_2, \]
 where $a_1,a_2,a_5$ are arbitrary constants. Then the Cauchy-Riemann equations (\ref{027}),
 (\ref{028}) imply that
 \[ \phi = - a_1x +a_5y + a_3, \]
 where $a_3$ is a constant. Now we substitute $\xi $, $\phi $ and $\tau $ into
 (\ref{032}) and (\ref{033}). The resulted equations, by a simple argument, imply
 \[\tau  =  2a_{3}x - 2a_{2}y + a_{4}+ 2a_{5}t. \]
 Finally, $\be =-2 a_5$ from (\ref{m09}). Summarizing, the infinitesimals are given by
 \bb\label{039} \left\{\begin{array}{r c l} \xi & = &
a_{1}y + a_{2}+a_{5}x, \\
\phi & = & - a_{1}x + a_{3} + a_{5}y,   \\
\tau & = & 2a_{3}x - 2a_{2}y + a_{4}+ 2a_{5}t, \\
\eta & = &- 2a_{5}, \end{array} \right.\ee which concludes the
proof of the main theorem in the case of exponential nonlinearity.

 \section{The Lie point symmetries for $f(u)=ku^{p}$}

 In this section we prove the main theorem in the case of nonlinearity of power type
 $f(u)=ku^p$. We suppose that $p\neq 0, p\neq 1, p\neq 2, p\neq 3$. We do not consider $p=2$ since
 in this case by a nonexistence result of Pokhozhaev and Veron \cite{pv} there is no solution of the
 corresponding Kohn-Laplace equation even in a very weak sense. The case $p=3$ will be treated in
 the next section. The case $p=1$ will be studied in section 9. Finally, if $p=0$ this is the item
 $(ii)$, which is reduced to $(i)$ as stated in the theorem.

 By (\ref{031}) we have
 \[ k(\al p +[2X\xi - \al ])u^p +kp\be u^{p-1} + (\lh \al )u+(\lh \be )=0.\]
 Hence $\be =0 $, $\lh\al =0$ and
 \bb\label{m14} \al =\frac{2}{1-p} X\xi . \ee
 By (\ref{m04}) and (\ref{m14}) it follows that
 \[ (p-3){\al }_t =0. \]
 Thus ${\al }_t =0$ since $p\neq 3$. Therefore $\al $ depends only on $x$ and $y$, and the
 equations (\ref{034}) and (\ref{035}) read
 \[ {\al }_x=-2{\phi }_t, \]
 \[ {\al }_y=2{\xi }_t . \]
 Hence there exist functions $B_1(x,y)$ and $B_2(x,y)$ such that
 \[ \xi =\frac{1}{2} {\al }_y (x,y) t +B_1(x,y), \]
 \[ \phi =-\frac{1}{2} {\al }_x (x,y) t +B_2(x,y). \]
 Substituting $\xi $ and $\phi $ into (\ref{027}) and (\ref{028}), we obtain
 \bb\label{m15} {\al }_{xy}=0 \ee
 and
 \bb\label{m16} {\al }_{xx}-{\al }_{yy}=0. \ee
 Since $\lh \al =0$ and $\al =\al (x,y)$, it follows that $\al $ is harmonic:
 \bb\label{m17} {\al }_{xx}+{\al }_{yy}=0. \ee
 From (\ref{m15}), (\ref{m16}) and (\ref{m17}) we conclude that
 \[ \al = Ax+By+C, \]
 where $A,B,C$ are constants. Thus
 \bb\label{m18} \xi =\frac{B}{2} t +B_1(x,y),\ee
 \bb\label{m19} \phi =-\frac{A}{2} t +B_2(x,y),\ee
 \bb\label{m20} \al =Ax+By+C. \ee
 Further, we substitute  $\xi $ and $\phi $ into (\ref{036}). In this way we see that ${\tau }_t$
 is a function of $x$ and $y$ only. Hence, there are functions $M(x,y)$ and $N(x,y)$ such that
 \bb\label{m21} \tau = M(x,y) t +N(x,y). \ee
 We put (\ref{m18}), (\ref{m19}) and (\ref{m21}) into (\ref{032}). We get that
 \[ M_x t +N_x + 2yM=g_1(x,y) - At, \]
 where $g_1$ is a function of $x$ and $y$ only. Thus $M_x = -A$. Hence
 \[ M=-Ax +m(y) \]
 for some function $m=m(y)$. Now we substitute (\ref{m18}) and (\ref{m21}) into (\ref{033}). We
 have
 \[ M_y t +N_y - 2xM=g_2(x,y) - Bt, \]
 where $g_2$ is a function of $x$ and $y$ only. Thus $M_y=-B$. Hence $m'(y)=-B$ and $m(y)=-By+D$,
 $D=$constant. Therefore
 \bb\label{m22} \tau = (-Ax-By +D) t +N(x,y). \ee
 By (\ref{m14}) and (\ref{m20}) we obtain
 \bb\label{m23} B_{1,x} = c_1x+c_2y +c_3 \ee
 where $c_1=(1-p)A/2$, $c_2=-(p+1)B/2$, $c_3=(1-p)/2$. By integration
 \bb\label{m24} B_1= \frac{c_1}{2} x^2 +c_2 xy +c_3x + \varphi (y) \ee
 for some function $\varphi $ of $y$ only. From (\ref{029}), (\ref{m18}), (\ref{m20})
 and (\ref{m24}) we obtain
 \[ {\varphi }''(y) = (3+p)A/2 =:c_4. \]
 Hence \[ \varphi (y) = c_4y^2 +c_5y +c_6, \]
 where $c_5,c_6$ are constants. Then
 \bb\label{m25}  B_1= \frac{c_1}{2} x^2 +c_2 xy +c_3x + c_4y^2 +c_5y +c_6. \ee
 By (\ref{027}), (\ref{m18}), (\ref{m19}) and (\ref{m25}):
 \[ B_{2,y}= d_1x +d_2y + c_3, \]
 where $d_1=-(p+1)A/2$, $d_2=(1-p)B/2$. Thus
 \bb\label{m26} B_2= d_1xy+\frac{d_2}{2}y^2 + c_3y +\psi (x) \ee
 for some function $\psi (x)$. From (\ref{030}), (\ref{m19}) and (\ref{m26}):
 \[ {\psi }''(x)=(3+p)B/2 =:c_7\]
 and therefore
 \[ \psi (x)=\frac{c_7}{2}x^2 +c_8x + c_9, \]
 where $c_8,c_9$ are constants. Hence
 \bb\label{m27} B_2= d_1xy+\frac{d_2}{2}y^2 + c_3y +\frac{c_7}{2}x^2 +c_8x + c_9. \ee
 Substituting (\ref{m18}) with $B_1$ given in (\ref{m25}) and
 (\ref{m19}) with $B_2$ given in (\ref{m27}) into (\ref{028}) we get that $c_8=-c_5$.

 It remains to determine the function $N(x,y)$ in (\ref{m22}). For this purpose we substitute $\tau $
 from (\ref{m22}) into (\ref{032}) and (\ref{033}), taking into account the
 already found expressions for $\xi ,\phi ,\al ,B_1$ and $B_2$. In this way we obtain
 \[ N_x= (2c_2 + 2c_7 -2B) x^2 +( 2 c_2 +d_2 +4 B)y^2 + (2c_1 +2 d_1 +2 A+4c_4) xy +
 (4c_3-2D)y + 2c_9, \]
 \[ N_y= (-5c_1 + 2d_1) x^2 +( 2A -2 c_4 -2d_1 )y^2 + (2 d_2 - 6 c_2-6 B -4c_7) xy +
 (2D-4c_3)x - 2c_6 .\]
 This system can be solved if and only if \[ D=2c_3\]
 and $A=B=0$. Hence \[ c_1=c_2=c_4=c_7=d_1=d_2=0. \] and
 the system is reduced to
 \[ N_x=2 c_9, \]
 \[ N_y=-2c_6, \]
 whose solution is $N=2c_{9}x -2 c_6 y +c_{10}$. After renaming the constants we obtain
 \bb\label{041} \left\{\begin{array}{r c l} \xi & = &a_{1}y + a_{2}+a_{5}x,  \\
\phi & = & - a_{1}x + a_{3} + a_{5}y,   \\
\tau & = & 2a_{3}x - 2a_{2}y + a_{4}+ 2a_{5}t, \\
\eta & = &\frac{2}{1-p} a_{5} u. \end{array} \right.\ee
 Observe that the dilation $Z$ comes from the constant $a_5$, while the rest corresponds to the
 generators in (\ref{u01}).

 \section{The Lie point symmetries for $f(u)=ku^{3}$}

 In this section we prove the second part of item $(iv)$ of the main theorem.

 Let $f(u)=ku^3$. Then $\be =0 $,
 \bb\label{m28} \al =-X\xi , \ee
 and
 \bb\label{m29} \lh\al =0 \ee
 from (\ref{031}). Applying $X$ to (\ref{034}), $Y$ to (\ref{035}) and adding, we obtain
 \[ \lh\al +2X{\phi }_t - 2  Y{\xi }_t=0.\]
 (Above we used (\ref{m01}) and (\ref{m06}).) By (\ref{m29})
 \[ X{\phi }_t -   Y{\xi }_t=0\]
 which together with (\ref{028}), differentiated with respect to $t$, implies
 \bb\label{m30} Y{\xi }_t=0 \ee
 and
 \bb\label{m31} X{\phi }_t=0. \ee
 Hence there exists a function $\var =\var (x,y)$ such that
 \bb\label{m32} Y\xi =\var \ee
 and, necessarily,
 \bb\label{m33} X\phi =-\var .\ee
 We also have by (\ref{m28}):
 \bb\label{m40} X\xi =-\al, \ee
 \bb\label{m41} Y\phi =-\al .\ee
 Then by (\ref{029}), (\ref{m40}) and (\ref{m41}):
 \[ 2X\al=\lh\xi =X^2\xi +Y^2\xi =X(-\al )+ Y\var = -X\al +{\var }_y, \]
 that is
 \[ 3X\al ={\var }_y.\]
 Hence, and from (\ref{034}), we obtain that there is a function $B_2(x,y)$ such that
 \bb\label{m42} \phi=-\frac{1}{6} {\var }_y\;t +B_2 (x,y). \ee
 Analogously
 \bb\label{m43} \xi=-\frac{1}{6} {\var }_x\;t +B_1 (x,y). \ee
 From (\ref{m43}) and (\ref{m30}) it follows that
 \bb\label{m44} {\var }_{xy}=0. \ee
 From (\ref{027}), (\ref{m43}) and (\ref{m42}) we have that
 \bb\label{m45} {\var }_{xx}-{\var }_{yy}=0. \ee
 Clearly, from (\ref{m44}) and (\ref{m45}), the function $\var $ is of the following form:
 \bb\label{m46} \var =\frac{k_1}{2} x^2 + \frac{k_1}{2} y^2 + k_2x +k_3y +k_4, \ee
 where $k_1,k_2,k_3,k_4$ are arbitrary constants.  In this way
 \bb\label{m47} \xi = -\frac{1}{6} (k_1x +k_2)\; t +B_1(x,y), \ee
 \bb\label{m48} \phi = -\frac{1}{6} (k_1y +k_3)\; t +B_2(x,y). \ee

 Substituting $\xi $ from (\ref{m47}) and $\varphi $ from (\ref{m46}) into (\ref{m32}), and
 integrating with respect to $y$, we obtain
 \bb\label{m49} B_1=\frac{1}{6}k_1x^2y + \frac{1}{6}k_1 y^3 + \frac{2}{3}k_2 xy +
 \frac{k_3}{2} y^2 +k_4 y +h_1(x), \ee
 where $h_1$ is a function of $x$ only. Analogously, from (\ref{m33}) we find
 \bb\label{m50} B_2=-\frac{1}{6}k_1x^3 - \frac{1}{6}k_1 xy^2 -\frac{k_2}{2} x^2
 - \frac{2}{3}k_3 xy  - k_4 x +h_2(y), \ee
 where $h_2$ is a function of $y$ only. After a substitution of $\xi $ and $\phi $ from
 (\ref{m47}) and (\ref{m48}) with $B_1$ and $B_2$ given by (\ref{m49}) and (\ref{m50}), into
 (\ref{027}), we obtain
 \[ {h_1}'(x)+\frac{1}{3}k_3x ={h_2}'(y)-\frac{1}{3}k_3y. \]
 Obviously, the last two equations can be easily integrated. In this way we find the functions $h_1$,
 $h_2$, and hence the functions $B_1$ and $B_2$. Summarizing, we have found
 \bb\label{m51}
 \begin{array}{lll}
 \xi & = & -\frac{1}{6} (k_1x +k_2)\; t +\frac{1}{6}k_1x^2y + \frac{1}{6}k_1 y^3 +
 \frac{2}{3}k_2 xy - \frac{1}{6}k_3x^2 +\frac{k_2}{2} y^2 \\
 & & +k_4 y +k_5x+k_6, \\
 \phi & = & -\frac{1}{6} (k_1y +k_3)\; t-\frac{1}{6}k_1x^3 - \frac{1}{6}k_1 xy^2 -\frac{k_2}{2} x^2
 - \frac{2}{3}k_3 xy  +\frac{1}{6} k_3y^2 \\
 & & - k_4 x+k_5y +k_7, \\
 \al & = & \frac{1}{6}k_1 t -\frac{1}{3}k_2y + \frac{1}{3}k_3 x - k_5.
 \end{array} \ee

 It remains to find $\tau $. In order to do this, we substitute (\ref{m51}) into
 (\ref{036}) and obtain
 \[ {\tau }_t = -\frac{1}{3}k_1 t + \frac{1}{3}k_2y -\frac{1}{3}k_3x +2k_5. \]
 Hence
 \bb\label{m52} \tau = -\frac{1}{6}k_1t^2 +(\frac{1}{3}k_2y -\frac{1}{3}k_3x +2k_5)t +N(x,y) \ee
 and the problem is reduced to the problem of finding the function $N$ in (\ref{m52}).
 Substituting (\ref{m51}) and (\ref{m52}) into equations (\ref{032}) and (\ref{033}), after some
 work, we finally obtain
 \[N_x= \frac{2}{3} k_1x^3 +\frac{2}{3}k_1xy^2 + k_2x^2 +\frac{2}{3}k_3 xy +
 \frac{1}{3}k_2y^2 + 2 k_7,  \]
 \[N_y= \frac{2}{3} k_1y^3 +\frac{2}{3}k_1x^2y + k_3y^2 +\frac{2}{3}k_2 xy +
 \frac{1}{3}k_3x^2 - 2 k_6.  \]
 The latter system can be easily solved. After renaming the constants, we have
\bb\label{043} \left\{\begin{array}{lll}
\xi & = & a_{1}(xt-x^{2}y-y^{3}) + a_{2}(t-4xy) +a_{3}(x^{2}-3y^{2}) +a_{4}x + a_{5}y  + a_{6}, \\
\phi & = & a_{1}(yt+x^{3}+xy^{2}) + a_{2}(3x^{2}-y^{2})  + a_{3}(t+4xy)+ a_{4}y  -a_{5}x  +a_{7},  \\
\tau & = & \displaystyle{a_{1}[t^{2}-(x^{2}+y^{2})^{2}]  +
a_{2}(-2yt-2x^{3}-2xy^{2})}\displaystyle{ + a_{3}(2xt-2x^{2}y-2y^{3})} \\
& & + 2a_{4}t  +2a_{7}x -2a_{6}y  +a_{8}, \\
\eta & = &-a_{1}t u +2a_{2}y u  -2 a_{3}x u -a_{4} u.
\end{array} \right.\ee
completing the proof of item $(iv)$ of the main theorem.

We observe that the dilation $Z$ is included in (\ref{043}).
Indeed, it corresponds to the constant $a_4$.

 \section{The Lie point symmetries for $f(u)=0$}

 The proof of item (i) of the main theorem is presented in this section. In order not to increase the
 volume of this paper, some of the calculations will be sketched, leaving the details to the
 interested reader.

 From (\ref{031}) with $f(u)=0$ we obtain
 \[ \lh \be =0 \]
 and
 \[\lh \al =0.\]
 From the latter equation we conclude, as in the beginning of section 7, that there exists a function
 $\var =\var (x,y)$ such that
 \bb\label{m60} Y\xi=\var , \ee
 and thus
 \bb\label{m61} X\phi =-\var \ee
 by (\ref{028}). On the other hand, from (\ref{m04}) it follows that there is a function $\psi =\psi
 (x,y)$ such that
 \bb\label{m62} \al = -X\xi +\psi , \ee
 \bb\label{m63} X\xi =-\al +\psi . \ee
 Following the arguments in obtaining (\ref{m42}) and (\ref{m43}) in the preceding section, we conclude
 that there exist functions $A=A(x,y)$ and $B=B(x,y)$ such that
 \bb\label{m64} \xi = (-{\var }_x+{\psi }_y)t/6 +A(x,y), \ee
 \bb\label{m65} \phi = -({\var }_y+{\psi }_x)t/6 +B(x,y). \ee
 Substituting $\xi $ and $\phi $ from (\ref{m64}) and (\ref{m65}) into (\ref{m60}), (\ref{m61}) and
 (\ref{027}), we obtain, respectively, that
 \bb\label{m66} {\psi }_{yy}={\var }_{xy}, \ee
 \bb\label{m67} {\psi }_{xx}=-{\var }_{xy}, \ee
 \bb\label{m68} {\var }_{yy}-{\var }_{xx} +2 {\psi }_{xy}=0. \ee
 Integrating (\ref{m66}) and (\ref{m67}) we have
 \bb\label{m69} {\psi }_y = {\var }_x + h_1(x), \ee
 \bb\label{m70} {\psi }_x=-{\var }_y +h_2(y) \ee
 for some functions $h_1=h_1(x)$ and $h_2=h_2(y)$. Then, differentiating (\ref{m69}) with respect to $x$
 and (\ref{m70}) with respect to $y$, adding and using (\ref{m68}), we get
 \[ h_1'(x)+h_2'(y)=0. \]
 Hence $h_1(x)=k_1x +k_2$ and $h_2(y)=-k_1y+k_3$ for some constants $k_1,k_2,k_3$. After renaming the
 constants and using (\ref{m69}) and (\ref{m70}), we obtain
 \bb\label{m71} \xi = (a_1x+a_2)t +A(x,y), \ee
 \bb\label{m72} \phi = (a_1x+a_3)t +B(x,y). \ee
 Further, from (\ref{m04}):
 \[ {\al }_t=-a_1. \]
 Hence
 \bb\label{m73} \al =-a_1t+g(x,y), \ee
 where the function $g=g(x,y)$ does not depend on $t$. Substituting
 (\ref{m71}), (\ref{m72}), (\ref{m73}) into (\ref{034}) and (\ref{035}) we find
 \bb\label{m74} \al =-a_1t-2a_3x+2a_2y+a_9, \ee
 where $a_9$ is an arbitrary constant and $a_1,a_2,a_3$ are the same which appear in (\ref{m71})
 and (\ref{m72}). Now, from (\ref{036}), (\ref{m71}) and (\ref{m72}), we deduce, after integration with
 respect to $t$, that there is a function $N=N(x,y)$ such that
 \bb\label{m75} \tau = a_1 t^2+(2A_x+4a_1xy+6a_2y-2a_3x)t +N(x,y). \ee
 We substitute $\xi $ from (\ref{m71}), $\phi $ from (\ref{m72}) into the determining
 equation (\ref{032}). In this way we obtain an identity which is linear in $t$. Equating to zero the
 corresponding coefficient of $t$, we obtain
 \bb\label{m76} A_{xx}=-2a_1y+2a_3 . \ee
 In an analogous way, using (\ref{033}),
 \bb\label{m77} A_{xy}=-2a_1x-4a_2. \ee
 We also have that (from equations (\ref{029}), (\ref{m71}) and (\ref{m72}))
 \bb\label{m78} A_{xx}+A_{yy}=-4a_3-8a_1y. \ee
 Then from (\ref{m76}), (\ref{m77}) and (\ref{m78}) we find
 \bb\label{m79} A=-a_1y^3-a_1x^2y+a_3x^2-3a_3y^2-4a_2xy+a_4x+a_5y+a_6, \ee
 where $a_5,a_6$ are constants.

 By substituting $\xi $ from (\ref{m71}) with $A$ given in (\ref{m79}), and $\phi $ from (\ref{m72}) into
 equations (\ref{027}) and (\ref{028}), and, then, integrating the resulted system for $B$, we find
 \bb\label{m80} B=a_1xy^2 +a_1x^3+4a_1xy +3a_2x^2-a_5x-a_2y^2+a_4y+a_7, \ee
 where $a_7$ is a constant.

 We have found $\xi ,\phi$ and $\eta $. To find $\tau $, it remains to determine the function $N(x,y)$ in
 (\ref{m75}).  From (\ref{032}), (\ref{033}), (\ref{m71}), (\ref{m72}), (\ref{m79}) and (\ref{m80}), we
 obtain the system
 \[ N_x = -4a_1x^3-4a_1xy^2 -6a_2x^2-4a_3xy-2a_2y^2+2a_7, \]
 \[ N_y = -4a_1x^2y-4a_1y^3-2a_3x^2-4a^2xy -6a_3y^2-2a_6, \]
 which can be easily solved. Our calculations can be summarized as
 \bb\label{044} \left\{\begin{array}{r c l} \xi & = &
a_{1}(xt-x^{2}y-y^{3})  + a_{2}(t-4xy) +a_{3}(x^{2}-3y^{2}) +a_{4}x + a_{5}y  + a_{6}, \\
\phi & = & a_{1}(yt+x^{3}+xy^{2})  + a_{2}(3x^{2}-y^{2})  + a_{3}(t+4xy)+ a_{4}y -a_{5}x +a_{7},  \\
\tau & = &a_{1}[t^{2}-(x^{2}+y^{2})^{2}]
+a_{2}(-2yt-2x^{3}-2xy^{2})
 + a_{3}(2xt-2x^{2}y-2y^{3})\\
 & & +  2a_{4}t +2a_{7}x -2a_{6}y +a_{8},\\
\eta & = &-a_{1}t u +2a_{2}y u -2 a_{3}x u +a_{9}u +\be(x,y,t).
\end{array} \right.\ee
where \bb\label{045} \lh\be=0, \ee $a_1,...,a_9$ are arbitrary
constants.

 \section{The Lie point symmetries for $f(u)=ku$}

 In this section we complete the proof of the main theorem.

 Let $f(u)=ku$, $k\neq 0$. Then by (\ref{031})
 \[ \lh \be +k\be =0\]
 and
 \bb\label{m81} \lh\al=-2kX\xi . \ee
 Applying $X$ to (\ref{034}) and $Y$ to (\ref{035}), and adding, we obtain
 \bb\label{m82} \lh \al =2Y{\xi }_t-2X{\phi }_t= 4Y{\xi }_t, \ee
 where we used (\ref{028}) and (\ref{m06}). Then from (\ref{m81}) and (\ref{m82}):
 \bb\label{m83} Y{\xi }_t=-kX\xi /2. \ee
 Hence and from (\ref{028}):
 \bb\label{m84} X{\phi }_t=kX\xi /2. \ee
 Differentiating (\ref{029}) we have
 \bb\label{m85} 2X{\al }_t=\lh{\xi }_t=X(X{\xi }_t)+Y(Y{\xi }_t)= X(-{\al }_t)-kY(X{\xi }_t)/2
 \ee
 by (\ref{m04}) and (\ref{m83}). On the other hand, by (\ref{m63})
 \bb\label{m86} X\xi = -\al +\psi ,\ee
 where $\psi =\psi (x,y)$. From (\ref{m85}) and (\ref{m86}) we have
 \[ 2X{\al }_t=-X{\al }_t +kY\al /2-k{\psi }_y /2 \]
 and hence
 \bb\label{m87} 3X{\al }_t = -k{\psi }_y /2 +kY\al /2. \ee
 Similarly
 \bb\label{m88} 3Y{\al }_t = k{\psi }_x /2 -kX\al /2. \ee
 We apply $X$ to (\ref{m87}), $Y$ to (\ref{m88}), and add:
 \[ 3(X^2+Y^2){\al }_t=k(XY-YX)\al /2 =k[X,Y]\al /2 = -2k{\al }_t, \]
 that is,
 \bb\label{m89} 3\lh {\al }_t=-2k{\al }_t.\ee
 Further, we differentiate (\ref{m81}) with respect to $t$ and use (\ref{m04}) to obtain
 \bb\label{m90} \lh{\al }_t =2k {\al }_t .\ee
 Since $k\neq 0$, from (\ref{m89}) and (\ref{m90}) it follows that ${\al }_t=0$ and hence $\al
 =\al (x,y)$. Thus, from (\ref{034}) and (\ref{035}), there exist functions $A=A(x,y)$ and
 $B=B(x,y)$ such that
 \bb\label{m91} \xi={\al }_y t/2 +A(x,y), \ee
 \bb\label{m92} \phi = -{\al }_x t/2 +B(x,y). \ee
 Since $\al $ does not depend on $t$, by (\ref{m91}) and (\ref{m04}) we have
 \bb\label{m93} {\al }_{xy}=0. \ee
 From (\ref{028}), (\ref{m91}) and (\ref{m92})
 \bb\label{m94} {\al }_{yy}-{\al }_{xx}=0. \ee
 The equations (\ref{m93}) and (\ref{m94}) can be easily solved. The solution is
 \bb\label{m95} \al = k_1x^2 +k_1y^2 +2k_3 x+2k_2y +k_4, \ee
 where $k_1,k_2,k_3,k_4$ are arbitrary constants. Hence
 \bb\label{m96} \xi=(k_1y+k_2) t +A(x,y), \ee
 \bb\label{m97} \phi = -(k_1 x+k_3) t +B(x,y). \ee
 On the other hand, from (\ref{036}), in which (\ref{m96}) and (\ref{m97}) are substituted,
 after an integration with respect to $t$, we obtain
 \bb\label{m98} \tau = [2y(k_1y+k_2) + 2x(k_1x+k_3) + 2(A_x+2k_1y^2+k_2y) ]t + N(x,y), \ee
 where the function $N=N(x,y)$ is to be determined.
 Further we substitute (\ref{m96}), (\ref{m97}) and (\ref{m98}) in (\ref{032}) and
 (\ref{033}). In this way we obtain two identities, linear in $t$. Equating the corresponding
 coefficients of $t$ implies
 \bb\label{m99} A_{xx}=-2k_1x -2 k_3, \ee
 \bb\label{m100} A_{xy}=-3k_2. \ee
 We observe now that equation (\ref{m81}) reads
 \bb\label{m101} 4k_1=-2kA_x-4k(k_1y^2+k_2y). \ee
 Differentiating (\ref{m101}) with respect to $x$ we obtain $0=-2kA_{xx}$ and hence
 \bb\label{m102} A_{xx}=0 \ee
 since $k\neq 0$. From (\ref{m99}) and (\ref{m102}) it follows that $k_1=k_3=0$.
 Then from (\ref{m101}), since $k\neq 0$, we have $A_x=-2k_2y$ which, together with
 (\ref{m100}) implies that $k_2=0$ and hence $A_x=0$. that is, $A=A(y)$. Summarizing we have
 obtained that
 \[ \xi = A(y),\;\;\;\phi =B(x,y),\;\;\;\tau =N(x,y),\;\;\;\al=k_4. \]
 Now the arguments in the section 4 imply that the infinitesimals are given by
 \bb\label{046} \left\{\begin{array}{r c l} \xi & = &
 a_{1}y + a_{2}, \\
 \phi & = & - a_{1}x + a_{3},  \\
 \tau & = & 2a_{3}x  - 2a_{2}y + a_{4},\\
 \eta & = & a_{5}u+\be(x,y,t),   \end{array} \right.\ee
 with
 \bb\label{047} \lh\be+k\be=0. \ee
 This completes the proof of the theorem.

\begin{center}\textbf{Acknowledgements} \end{center}We thank Enzo Mitidieri for his suggestion that we write this
 paper as well as
 for his firm encouragement. Y. Bozhkov would also like to thank
 FAPESP and CNPq, Brasil, for financial support. I. L. Freire is grateful to
 CAPES, Brasil, for financial support.

\label{lastpage}

\end{document}